\documentclass[12pt,draft]{article}
\usepackage{amsmath,amsfonts,calrsfs,amssymb, color}
\newtheorem{Theorem}{Theorem}[section]

\newtheorem{Proposition}[Theorem]{Proposition}
\newtheorem{Lemma}[Theorem]{Lemma}

\newtheorem{Remark}[Theorem]{Remark}

\newtheorem{Hypothesis}[Theorem]{Hypothesis}
\makeatletter
\newif\ifmsbmloaded@

\@addtoreset{equation}{section}

\makeatother

\newcommand{\one}{1\!\!\!\;\mathrm{l}}

\def\R{\mathbb R}
\def\N{\mathbb N}

\def\E{\mathbb E}
\def\P{\mathbb P}
\def\Q{\mathbb Q}
\def\ds{\displaystyle}

\title{\bf Existence and uniqueness of solutions for Fokker--Planck equations on Hilbert spaces}
\author{Vladimir Bogachev\thanks{Supported in part by
the RFBR project
 07-01-00536,
the Rus\-sian--Japanese Grant 08-01-91205-JF, the
Rus\-sian--Ukrainian Grant 08-01-90431,
SFB 701 at the University of Bielefeld},\\
Department of Mechanics and Mathematics,\\
          Moscow State University, 119991 Moscow, Russia,\\\\
Giuseppe Da Prato,\\
 Scuola Normale Superiore
di Pisa, Italy\\
and\\
 \\
Michael R\"ockner \thanks{Supported by the DFG through SFB-701 and IRTG 1132 as well as the
BIBOS-Research Center. } \\
Faculty of Mathematics, University of Bielefeld, Germany\\ and\\
Department of Mathematics and Statistics,\\ Purdue University, W. Lafayette, 47906, IN,   U. S. A.
}
\date{June 12 2009}

\begin{document}

\maketitle
\begin{abstract}
We consider a stochastic differential equation in a Hilbert space with time-dependent coefficients for which no general existence and uniqueness results are known. We prove, under suitable assumptions, existence and uniqueness of a measure valued solution, for the corresponding Fokker--Planck equation. In particular, we verify the Chapman--Kolmogorov equations and get an evolution system of  transition probabilities for the stochastic dynamics informally given by the stochastic differential equation.

 \end{abstract}

\noindent {\bf 2000 Mathematics Subject Classification AMS}: 60H15, 60J35, 60J60, 47D07


\bigskip

\maketitle

\section{Introduction}
In recent years there has been a lot of interest  in Fokker--Planck and transport equations with
irregular coefficients in finite dimensions (see e.g.
\cite{A04}, \cite{A05}, \cite{L07}, \cite{M07},
\cite{Fi08}, \cite{LeLi08}  and the references therein  and also the fundamental paper \cite{DPL89}).
More recently, also transport equations in infinite dimensions  have been analyzed
(see, e.g., \cite{0}, \cite{BDPRS09}). In \cite{BDPR08c}, \cite{BDPR09}
we have started a study of Fokker--Planck equations
in infinite dimensions, more precisely, on Hilbert spaces. In the present paper we continue
this study by proving existence and uniqueness results for irregular (even non continuous)
drift coefficients.
Here we consider the case of full noise (i.e. the diffusion operator is invertible).
Another paper concerned with
degenerate (Hilbert--Schmidt) noise is in preparation. The case of zero noise, even when the drift
coefficients depends (nonlinearly) on the solutions is treated in finite dimensions in \cite{BRS09}
and in infinite dimensions in \cite{BDPRS09}.\medskip

Before we describe our framework and results more precisely, we would like to stress
that we can also prove the Chapman--Kolmogorov equations for our solutions
to the Fokker--Planck equations. This is, of course, a consequence of uniqueness
of solutions, which in turn follows from a technique developed by us in several
 papers first in finite (see \cite{BDPRS07} and also \cite{BRS00},
 \cite{BRS02} for the elliptic case) and subsequently in infinite
 dimensions  (see \cite{BDPR08c} and Section 3 below).\bigskip

Let $H$ be a separable real Hilbert space with
inner product $\langle\cdot,\cdot   \rangle$ and corresponding norm $|\cdot|$.  $L(H)$ denotes the set of all bounded linear  operators  on $H$ with its usual norm $\|\cdot\|$, $ \mathcal B(H)$ its Borel $\sigma$-algebra, $ \mathcal \mathcal B_b(H)$ the set of all bounded $ \mathcal B(H)$-measurable functions
from $H$ to $\R$  and $ \mathcal P(H)$ the set of all probability measures on $H$, more precisely on $(H, \mathcal B(H))$.

Consider the following type of non-autonomous stochastic differential equation on $H$ and time interval $[0,T]$:
\begin{equation}
\label{e1.1}
\left\{\begin{array}{l}
dX(t)=(AX(t)+F(t,X(t)))dt+\sqrt{C}dW(t),\\
\\
X(s)=x\in H,\;t\ge s.
\end{array}\right.
\end{equation}
Here  $W(t),\;t\ge 0,$ is a cylindrical Wiener process on $H$ defined on a stochastic basis $(\Omega, \mathcal F, (\mathcal F_t)_{t\ge 0},\P)$, $C$ is a symmetric positive operator in $L(H)$,
$D(F)\subset\mathcal B([0,T]\times H),$
$F\colon D(F)\subset [0,T]\times H\to H$, $t\in [0,T],$ is a measurable map, and  $A:D(A)\subset H\to H$ is the infinitesimal generator of a $C_0$-semigroup $e^{tA},\;t\ge 0,$ in $H.$

Without further regularity assumptions on $F$ it is, of course, not at all clear whether \eqref{e1.1} has a solution in the strong or even in the weak sense.
If, however,
there is a weak solution to \eqref{e1.1}, then it is a well known consequence of It\^o's formula that its transition probabilities $p_{s,t}(x,dy),\;x\in H,\;s\le t$, solve the Fokker--Planck equation determined by  the associated Kolmogorov operator. The purpose of this paper is to describe very general conditions on $F$ above for which one can solve the Fokker--Planck equation directly for Dirac initial conditions and thus to obtain the transition functions $p_{s,t},\;s\le t$, corresponding to \eqref{e1.1} though one might not have a solution to it.
In particular, we prove that $p_{s,t},\;s\le t$, satisfy the Chapman--Kolmogorov equation under wide conditions.

The general motivation to study Fokker--Planck equations instead of Kolmogorov equations, as done in our previous papers, is that the latter are equations for functions, whereas the first are equations for measures for which one has e.g. much better compactness criteria in our infinite dimensional situation. So, there is a good chance to obtain very general existence results.\medskip

Before we write down the Fokker--Planck equation we
recall that the Kolmogorov operator $L_0$ corresponding to \eqref{e1.1} reads as follows:
\begin{multline*}
L_0u(t,x)=D_tu(t,x)+\frac12\;\mbox{\rm Tr}\;[CD_x^2u(t,x)]\\+\langle x,A^*D_xu(t,x)
 \rangle+\langle F(t,x),D_xu(t,x)  \rangle,\quad x\in H,\;t\in[0,T],
\end{multline*}
where $D_t$ denotes the derivative in time and $D_x,D^2_x$ denote the first and second order
Fr\'echet derivatives in space, i.e. in $x\in H$,  respectively.
The operator $L_0$ is defined  on the space $D(L_0):={\mathcal E}_{A}([0,T]\times H)$,   the linear span of all
 real   parts of    functions $u_{\phi,h}$ of the form
 \begin{equation}
\label{e1.3}
u_{\phi,h}(t,x)=\phi(t)e^{i\langle x,h(t) \rangle},\quad t\in [0,T], \;x\in H,
\end{equation}
where $\phi\in C^1([0,T]),\;\phi(T)=0$, $h\in C^1([0,T];D(A^*))$  and $A^*$ denotes the adjoint of $A$.

For a fixed initial time $s\in[0,T]$ the  Fokker--Planck equation   is an equation for measures  $\mu(dt,dx)$ on $[s,T]\times H$ of the type
$$
\mu(dt,dx)=\mu_t(dx)dt,
$$
with $\mu_t\in \mathcal P(H)$ for all $t\in[s,T],$ and $t\mapsto \mu_t(A)$ measurable on $[s,T]$ for all $A\in \mathcal B(H)$,
i.e., $\mu_t(dx),\;t\in[s,T],$ is a probability kernel from
$([s,T],\mathcal B([s,T])$ to $(H,\mathcal B  (H))$. Then the equation
for an initial condition $\zeta\in  \mathcal P(H)$ reads as follows:
$\forall\;u\in D(L_0)$ one has
\begin{multline}
\label{e1.5}
\int_Hu(t,y)\mu_t(dy)=\int_Hu(s,y)\zeta(dy)+\int_s^tds' \int_HL_0u(s',y)\mu_{s'}(dy)
,\\\quad \mbox{\rm for   $dt$-a.e.}\;t\in [s,T],
\end{multline}
where a $dt$-zero set may depend on $u$.
When writing \eqref{e1.5} (or \eqref{(1.6a)},
\eqref{(1.6b)} or \eqref{(1.8)} below) we always implicitly assume that
\begin{equation}
\label{e1.5'}
\int_{[0,T]\times H}(|\langle y,A^*h(t) \rangle|+|F(t,y)|)\mu(dt,dy)<\infty
\end{equation}
for all $h\in C^1([0,T];D(A^*))$, so that all involved integrals exist in the usual sense.

\begin{Remark}[Equivalent formulations]
\label{r1.1}
\em

\noindent (i) We would like to emphasize that a priori we do not assume any continuity of the map
$$
t\mapsto\int_H\varphi(y)\mu_t(dy),\quad t\in [s,T],
$$
for ``sufficiently many'' nice functions $\varphi:H\to \R$, as e.g. $\varphi\in \mathcal E_A(H),$ defined to be the set of linear combinations of all
 real   parts of    functions  of the form
$$
H\ni x\mapsto e^{i\langle x,h  \rangle},\quad h\in D(A^*).
$$
Nevertheless, one can prove that, under the assumption \eqref{e1.5'},
identity \eqref{e1.5}  is equivalent to the usual ``differential form'' of the Fokker--Planck equation:
$\forall\;u\in D(L_0)$, $\forall\;\varphi\in\mathcal E_A(H)$ one has
\begin{equation}
\label{(1.6a)}
\frac{d}{dt}\;\int_Hu(t,y)\mu_t(dy)=\int_HL_0u(t,y)\mu_t(dy),
\quad\mbox{\rm for $dt$-a.s. $t\in[s,t],$}
\end{equation}
\begin{equation}
\label{(1.6b)}
\lim_{t\to s}\int_H\varphi(y)\mu_t(dy)=
\int_H\varphi(y)\zeta(dy).
\end{equation}

Here (since no continuity is assumed on $t\mapsto\mu_t,\;t\in[s,T]$) the limit in \eqref{(1.6b)} has to be understood in the following sense: there exists a map
$t\mapsto\widetilde{\mu_t}\in\mathcal P(H),\;t\in[s,T],$ equal to $t\mapsto\mu_t$ outside a set of $dt$-measure zero so that \eqref{(1.6b)} holds with
$\widetilde{\mu_t}$ in place of $\mu_t$. That \eqref{e1.5} and \eqref{(1.6a)}+\eqref{(1.6b)}   are indeed equivalent, was proved in
\cite[Remark 1.2]{BDPR09}. Considering $D(L_0)$ as test functions and dualizing we then turn  \eqref{(1.6a)}+\eqref{(1.6b)} into the familiar form of the Fokker--Planck equation
$$
\frac{\partial}{\partial t}\;\mu_t=-L_0^*\mu_t,\quad \mu_s=\zeta.
$$

\noindent(ii) Setting $t=T$ and recalling that $u(T,\cdot)\equiv 0$ for all $u\in D(L_0)$ we see that
(under assumption \eqref{e1.5'}) equation \eqref{e1.5} is obviously also equivalent to
\begin{equation}\label{(1.8)}
\int_{[s,T]\times H}L_0 u(s',y)\mu(ds',dy)=-\int_Hu(s,y)\zeta(dy),\quad\forall\;u\in D(L_0).
\end{equation}

\noindent(iii) By an easy approximation argument it follows that if \eqref{e1.5} holds for all $u\in D(L_0)$, then it holds for all $u$ of the form \eqref{e1.3} with $h\in C([0,T];D(A^*))$ and $h=h_1+\cdots+h_N$ with $h_i\in C^1([s_{i-1},s_i];D(A^*)),\;1\le i\le N$, and $0=s_0<s_1<\cdots<s_N=T.$
\end{Remark}
Solving \eqref{e1.5} (if this is possible) with $\zeta=\delta_x$ (:=Dirac measure in $x\in H$) for $x\in H$ and $s\in [0,T)$ and expressing the dependence on $x,s$ in the notation, we obtain probability measures $p_{s,t}(x,dy),\;t\in[s,T],$
such that the measure $p_{s,t}(x,dy)dt$ on $[s,T]\times H$ is a solution of \eqref{e1.5}. We shall see in Section 3 below, that if we have
uniqueness for \eqref{e1.5} and a ``sufficient continuity'' of the functions
 $t\mapsto p_{s,t}(x,dy)$, then these measures satisfy the Chapman--Kolmogorov equations, i.e. for $0\le r<s<t\le T$ and $x\in H$ (or in a properly chosen subset thereof)
\begin{equation}
\label{(1.9)}
\int_Hp_{s,t}(x',dy)p_{r,s}(x,dx')=p_{r,t}(x,dy),
\end{equation}
where the left hand side is a measure defined for $A\in \mathcal B(H)$ as
$$
\int_{H\times H}\one_A(y)p_{s,t}(x',dy)p_{r,s}(x,dx').
$$
The theoretical component of the  paper consists of two parts. In the first part (see Section 2 below) we shall prove existence of solutions to \eqref{e1.5} under very general assumptions on coefficients $A, F$ and $C$. There is
a well known generic difference between the case
when $C$ has finite trace or not. We shall concentrate on the latter, more precisely, even on the
extreme situation when $C^{-1}\in L(H)$ (hence including the ``white noise'' case).
The reason is that if Tr $C<\infty$, there are a number of known
existence  results  (cf. \cite{BDPR08b} and also \cite{BDPR04}, \cite{BDPR08})
based on the method of constructing Lyapunov functions with weakly compact level sets
for the Kolmogorov operator $L_0$, which does not apply when Tr $C=\infty$.
We refer to Theorem 2.5 below for the precise formulation of our result and to Remark 2.3(ii) for the relations of our method with Lyapunov functions. We only emphasize here that under the assumptions of Theorem 2.5,
on the one hand we are very far away from being allowed to apply Girsanov--Maruyama's theorem to weakly solve \eqref{e1.1}, whereas, on the other hand, the proof of  Theorem 2.5  heavily relies on  applying Girsanov--Maruyama's transformation to a proper approximation. Furthermore,
we only need the continuity of the components $x\mapsto \langle h,F(x)\rangle,\;h\in H,$ of $F$ (but see also Remark 2.6(ii) below).\medskip

The second part of the paper (see Section 3) is devoted to uniqueness of solutions to \eqref{e1.5} and to
deriving the Chapman--Kolmogorov equations \eqref{(1.9)}. Here additional dissipativity (not continuity) conditions on $F$ are needed and we rely heavily on the uniqueness results in
\cite{BDPR08c}, which hold no matter whether $C$ is of trace class or not, hence also apply in case of the existence results of \cite{BDPR08b}.\medskip

In the last part (see Section 4) we present applications which, in particular include reaction-diffusion equations with polynomially growing, time dependent nonlinearities.\medskip

Finally, we would like to mention that some of our results in Section 2 have been announced (though in a   weaker formulation) in \cite{BDPR09} with rough sketches of the proofs.

\section{Existence of solutions of the Fokker--Plank equation}

Let us first introduce some assumptions to be used below.

\begin{Hypothesis}
\label{h2.1}
\begin{enumerate}
\item[]

\item[(i)] $A$ is  self-adjoint and such that there exists $\omega\in \R$ such that
$$\langle Ax,x  \rangle\le \omega|x|^2,\;x\in D(A).$$

\item[(ii)] $C\in L(H)$ is   symmetric,  nonnegative and such that $C^{-1}\in L(H)$.

\item[(iii)]   There exists $\delta\in (0,1/2)$ such that $(-A)^{-2\delta}$ is of trace class.
 \end{enumerate}
 \end{Hypothesis}

 Let us notice that  it follows from (iii) that
 the embedding $D(A)\subset H$ is compact.

 It is well known that, under Hypothesis \ref{h2.1},  the stochastic convolution
$$
W_A(t)= \int_0^te^{(t-s)A}\sqrt C dW(s),\quad t\ge 0,
$$
is a well defined mean square continuous process in $H$ with values in $D((-A)^{\delta})$ and that
\begin{equation}
\label{e2.1}
\sup_{t\in[0,T]}\E|(-A)^{\delta}W_A(t)|^2\le\|C\|\;\mbox{\rm Tr}\;[(-A)^{-2\delta}]:=c_\delta.
\end{equation}

\begin{Hypothesis}
\label{h2.1r}
There exist bounded measurable maps \mbox{$F_\alpha\colon\, [0,T]\times H\to H$,} $\alpha\in(0,1],$  such that for all $(t,x)\in D(F)$ and all $h\in D(A)$
$$
\lim_{\alpha\to 0}\langle h,F_\alpha(t,x)\rangle=\langle h, F(t,x)\rangle,
$$
\begin{equation}\label{e2.4}
|F_\alpha(t,x)|\le |F(t,x)|,
\end{equation}
\begin{equation}\label{e2.5}
|\langle h,F(t,x)-F_\alpha(t,x)\rangle|\le \alpha c(h)|F(t,x)|,
\end{equation}
for some constant $c(h)>0$.
\end{Hypothesis}
Now we consider the following approximating stochastic equations for fixed $s\in [0,T]$:
\begin{equation}\label{e2.6}
\left\{\begin{array}{l}
dX_\alpha(t)=[AX_\alpha(t)+F_\alpha(t,X_\alpha(t))]dt+\sqrt{C}dW(t),\\
\\
X_\alpha(s)=x,\quad s\le t.
\end{array}\right.
\end{equation}
 Since $C^{-1}\in L(H)$, by Girsanov's theorem it follows that for every $x\in H$ equation \eqref{e2.6}
 has a martingale solution which we denote by $X_\alpha(\cdot,s,x)$
 (see e.g. \cite[Proposition 10.22]{DPZ1}). Let   $W(t),\;t\ge s,$ denote the corresponding
 cylindric Wiener process on $H$ and set
$$
W_A(t,s)=\int_s^te^{(t-s')A}\sqrt C\;dW(s'),\quad t\ge s.
$$

Let us introduce the transition evolution operator
$$
P_{s,t}^\alpha\varphi(x)=\E[\varphi(X_\alpha(t,s,x))],\quad 0\le s<t\le T,\;\varphi\in \mathcal B_b(H).
$$
The    Kolmogorov operator $L_\alpha$ corresponding to \eqref{e2.6} is  given by
the following expression for $u\in D(L_0)$:
\begin{multline*}
L_\alpha u(t,x)=D_tu(t,x)+\frac12\;\mbox{\rm Tr}\;[CD_x^2u(t,x)]\\+\langle x,A^*D_xu(t,x)  \rangle+\langle F_\alpha(t,x),D_xu(t,x)  \rangle,\quad x\in H,\;t\in [0,T].
\end{multline*}
From now on we fix $s\in[0,T)$ and set
$$
\mu_t^\alpha(dx):=(P_{s,t}^\alpha)^*\zeta(dx),
$$
where $\zeta\in\mathcal P(H)$ is the initial condition. So,
$$
\int_H\varphi(y)\mu_t^\alpha(dy)=\int_HP_{s,t}^\alpha\varphi(y)\zeta(dy),\quad\forall\;\varphi\in \mathcal B_b(H).
$$

Then by It\^o's formula  this gives a solution to the corresponding Fokker--Planck equation
\begin{multline}\label{e2.9}
\int_Hu(t,x)\mu^\alpha_t(dx)=\int_Hu(s,x)\zeta(dx)+\int_s^tds \int_HL_\alpha u(s',x)\mu^\alpha_{s'}(dx)
,\\\quad \mbox{\rm for all}\;t\in [s,T],\;\forall\;u\in D(L_0).
\end{multline}
Now we introduce our crucial assumption.

\begin{Hypothesis}
\label{h2.2}
There exist $K>0$ and a lower semicontinuous function
$V\colon\, [s,t]\times H\to [1,\infty]$ such that $|F|\le V$ on $[s,T]\times H$, where here and below $|F|=\infty$ on $([s,T]\times H)\setminus D(F)$, and
\begin{equation}\label{e2.10}
P^\alpha_{s,t}V^2(t,\cdot)(x)\le KV^2(t,x)<\infty,\quad \forall\;(t,x)\in D(F),\;t\in [s,T],\;\alpha\in (0,1]
\end{equation}
\end{Hypothesis}

\begin{Remark}
\em (i) Since we can always add a constant to $V$ preserving all its properties, the assumption that $V\ge 1$ is not a restriction. Furthermore, \eqref{e2.10} implies that
\begin{equation}\label{e2.10'}
P^\alpha_{s,t}\one_{H\setminus D(F(t,\cdot))}(x)=0\quad\forall\;  (t,x)\in D(F),\;t\in [s,T],\;\alpha\in(0,1],
\end{equation}
where
$$
 D(F(t,\cdot))=\{x\in H:\;\exists\;t\in[0,T]\;\mbox{\rm such that}\;(t,x)\in D(F)\}.
$$

(ii) Roughly speaking to satisfy Hypothesis \ref{h2.2} means that we have to find a function which is a Lyapunov function for $P^\alpha_{s,t}$ (not for $L_0$ as in \cite{BDPR08b}) uniformly in $\alpha$, and whose square root dominates the nonlinear part of the drift of \eqref{e1.1}.
\end{Remark}

\begin{Lemma}\label{l2.4}
 Assume that  Hypotheses \ref{h2.1}
and \ref{h2.2} hold. Then for all $\alpha\in(0,1],$ $\zeta\in\mathcal P(H)$, $t_1,t_2\in[s,T]$
one has
$$
\int_{t_1}^{t_2}\int_H V^2(s',x) \mu^\alpha_{s'}(dx)ds'\le K\int_{t_1}^{t_2}\int_HV^2(s',x) \zeta(dx)ds'.
$$
In particular, if
$$
\int_s^T\int_HV^2(s',x) \zeta(dx)ds'<\infty,
$$
then
$$
\int_s^T\int_HP^\alpha_{s,t}\one_{H\setminus D(F(t,\cdot))}(s',x)\zeta(dx)dt=0, \quad\forall\; \alpha\in(0,1].
$$
\end{Lemma}
{\bf Proof}. The first assertion is an immediate consequence of \eqref{e2.10}.
The second then follows from \eqref{e2.10'} since $V=\infty$ on $([s,t]\times H)\setminus D(F)$.
Hence
$$
\int_s^T\int_H  \one_{([s,T]\times H)\setminus D(F)}(s',x)\zeta(dx)ds'=0
$$
by our assumption. $\Box$

Now we can state and prove our main existence result.
 \begin{Theorem}
\label{t2.5}
 Assume that  Hypotheses \ref{h2.1}-\ref{h2.2} hold and that
\begin{multline}\label{e2.11}
(t,x)\mapsto \langle h,F^\alpha(t,x)\rangle \;\mbox{\it is continuous on}\; [s,T]\times H,\;\;\forall\; h\in D(A),\; \alpha\in(0,1].
\end{multline}
Let $\zeta\in \mathcal P(H)$ be such that
\begin{equation}
\label{e2.12}
\int_s^T\int_H(V^2(s',x)+|x|^2)\zeta(dx)ds'<\infty.
\end{equation}
Then there exists a solution $\mu_t(dx)dt$ to the Fokker-Planck equation \eqref{e1.5}
such that
$$
\sup_{t\in[s,T]}\int_H|x|^2\mu_t(dx)<\infty
$$  and
 $$t\mapsto \int_Hu(t,x)\mu_t(dx)$$ is continuous
on $[s,T]$ for all $u\in D(L_0)$.
 In particular, \eqref{e1.5} holds for all $t\in[s,T]$. Finally, for some $C>0$ one has
\begin{equation}
\label{e2.12'}
\begin{array}{l}
\ds\int_s^T\int_H \bigl(V^2(s',x)+|(-A)^\delta x|^2+|x|^2\bigr)\, \mu_{s'}(dx)ds'\\
\\
\ds\le C
\int_s^T\int_H \bigl(V^2(s',x)+|x|^2\bigr)\, \zeta(dx)ds'
\end{array}
\end{equation}
and hence $\mu_t(D(F(t,\cdot)))=1$ for all $dt$-a.e. $t\in[s,T]$.
\end{Theorem}

\begin{Remark}
\label{r2.6}
\em (i) The idea to prove the above result is to show that
the measures
$\mu^\alpha_t(dx)dt$, $\alpha\in(0,1]$, on $[0,T]\times H$ are uniformly
tight and that a limit point solves  \eqref{e1.5}. Only for the latter part (i.e. Claim 3 of the proof of Theorem \ref{t2.5} below) condition \eqref{e2.11} is needed.\medskip

\noindent (ii) We believe that, in fact, condition \eqref{e2.11} is superfluous. This was proved in \cite[Theorem 5.2]{6} in the time independent case.
Some of the ingredients of the proof have, however, not yet been proved in the time dependent case though they are very likely to hold also here. A corresponding paper is  in preparation.\medskip

\noindent (iii) If the continuity condition \eqref{e2.11} can be dropped, then so can Hypothesis \ref{h2.1r}. Simply define for $\alpha\in(0,1]$
$$
F_\alpha(t,x)=\left\{\begin{array}{l}
\frac{F(t,x)}{1+\alpha|F(t,x)|}\quad\mbox{\rm if $(t,x)\in D(F)$},\\\\
0\quad\mbox{\rm otherwise}.
\end{array}\right.
$$
Then obviously $F_\alpha$ enjoys all properties in Hypothesis \ref{h2.1r}.
\end{Remark}

\noindent{\bf Proof of Theorem \ref{t2.5}}. Below we shall use
the weak topology $\tau_w$ on $H$ and weak convergence of a sequence
of measures $\nu_n\in\mathcal P(H)$ on $(H,|\cdot|)$ and on $(H,\tau_w)$.
To avoid confusion we shall use the terminology
``weak convergence of $\nu_n$'' as usual if we refer to the norm topology of $H$ and
``$\tau_w$-weak convergence of $\nu_n$''
if we refer to $\tau_w$. Here we recall that since $H$ is always assumed to be separable, the Borel $\sigma$-algebra with respect to $\tau_w$ coincides with $\mathcal B(H)$.

The proof is structured in three claims.\medskip

\noindent{\bf Claim 1}. For any given sequence in $(0,1]$
convergent to zero
there exists a subsequence $\alpha_n\to 0$   and measures $\mu_t,\;t\in[0,T],$ such that the measures
$\mu_t^{\alpha_n}$ converge
$\tau_w$-weakly  to $\mu_t$ for all $t\in[0,T]$. Furthermore,
$$\sup_{t\in[s,T]}\int_H|x|^2\mu_t(dx)<\infty$$ and  for all $u\in D(L_0)$ the map
$$t\mapsto\int_Hu(t,x)\mu_t(dx)$$ is continuous on $[s,T]$. In particular,
$\mu_t(dx),\;t\in[s,T],$ are probability kernels from $([s,T],\mathcal B([s,T]))$
to $(H,\mathcal B(H)$.\medskip

\noindent{\bf Claim 2}. Selecting another subsequence we may assume that
the measures $\mu_t^{\alpha_n}(dx)dt$ converge weakly to $\mu_t(dx)dt$ on $[0,T]\times H$
where $\mu_t(dx),\;t\in [0,T],$ is defined as in Claim~1.
Furthermore, \eqref{e2.12'} holds.\medskip

\noindent{\bf Claim 3}. The measure
$\mu_t(dx)dt$ from Claim 2 solves the Fokker--Planck equation \eqref{e1.5}.
\medskip

\noindent{\bf Proof of Claim 1}. Let $\alpha\in(0,1]$, set $X_\alpha(t):=X_\alpha(t,s,x),\;x\in H,$ and
\begin{equation}
\label{e2.13}
Y_\alpha(t):=X_\alpha(t)-W_A(t,s),\quad t\ge s.
\end{equation}
Then in the mild sense
$$
\frac{d}{dt}\;Y_\alpha(t)=AY_\alpha(t)+F_\alpha(t,X_\alpha(t)),\quad t>s.
$$
Multiplying both sides by $Y_\alpha(t)$ for $t>s$ we obtain
$$
\frac12\;\frac{d}{dt}\;|Y_\alpha(t)|^2+|(-A)^{1/2}Y_\alpha(t)|^2= \langle F_\alpha(t,X_\alpha(t)),Y_\alpha(t)   \rangle.
$$
Integrating over $[s,T]$ and applying Young's inequality we get that for $t\ge s$
\begin{multline}
\label{e2.14}
|Y_\alpha(t)|^2+2\int_s^t|(-A)^{1/2}Y_\alpha(s')|^2ds'
\\
\le |x|^2+\int_s^t(|Y_\alpha(s')|^2+|F_\alpha(s',X_\alpha(s'))|^2)ds'.
\end{multline}
The above derivation of \eqref{e2.14} is a bit informal since $A$ is in general unbounded. This can, however,  easily be made rigorous by approximation (see \cite[Section 3.27]{DP2004}). Dropping the term involving $A$ and applying Gronwall's lemma we deduce from \eqref{e2.14} that for $t\ge s$
\begin{equation}
\label{e2.14'}
|Y_\alpha(t)|^2\le e^{t-s}|x|^2+\int_s^te^{t-s'}|F_\alpha(s',X_\alpha(s'))|^2ds'.
\end{equation}
Taking expectation and applying \eqref{e2.4} and Hypothesis \ref{h2.2}, yields
$$
\E|Y_\alpha(t)|^2\le e^{t-s}|x|^2+K\int_s^te^{t-s'}|V^2(s',x)|^2ds',\quad t\ge s
$$
and after resubstituting according to \eqref{e2.13} it follows that for $s\le t\le T$
$$
\E|X_\alpha(t,s,x)|^2\le 2e^{T-s}|x|^2+2Ke^{T-s}\int_s^T|V^2(s',x)|^2ds'+2\kappa,
$$
where
$$
\kappa:=\sup_{t\in[s,T]}\E|W_A(t)|^2 (<\infty).
$$
Now we integrate with respect to $\zeta$ over $x\in H$ and obtain for $s\le t\le T$
\begin{equation}
\label{e2.15}
\int_H|x|^2\mu^\alpha_t(dx)\le C \left[1+\int_s^T\int_H(V^2(s',x)+|x|^2)   ds'\zeta(dx)\right],
\end{equation}
for some $C>0$. By \eqref{e2.12} the right hand side of \eqref{e2.15} is finite. But it is also independent of $\alpha\in(0,1]$. Consequently, since closed balls in $H$ are $\tau_w$-compact and metrizable we can apply a version of Prohorov's theorem on completely regular topological spaces (see \cite[Theorem 8.6.7]{B2007})
which implies that given any  sequence in $(0,1]$ convergent to zero,
for each $t\in[s,T]$, there exists a sub-sequence $\{\alpha_n\}$ (dependent on $t$)
such that the measures $\mu_t^{\alpha_n}$ converges $\tau_w$-weakly
to a measure $\widetilde{\mu_t}\in\mathcal P(H)$ as $n\to \infty$.

To prove that this sequence $\{\alpha_n\}$ can indeed be chosen independently of $t\in[s,T]$ we need to prove that for each $\varphi\in \mathcal E_A(H)$ and
$$
\mu_t^\alpha(\varphi):=\int_H\varphi(x)\mu_t^\alpha(dx),\quad t\in[s,T],\;\alpha\in(0,1]
$$
we have:
\begin{equation}
\label{e2.17}
\mbox{\rm the maps}\;t\mapsto\mu_t^\alpha(\varphi),\;\alpha\in(0,1], \mbox{\rm are equicontinuous on}\;[s,T].
\end{equation}
Suppose \eqref{e2.17} is true for all $\varphi\in \mathcal E_A(H)$, then we can proceed as follows.
By a diagonal argument we can choose $\{\alpha_n\}$ such that $\mu_t^{\alpha_n}\to \widetilde{\mu_t}$
 $\tau_w$-weakly as $n\to \infty$ for every rational $t\in[s,T]$. We note that since $|\cdot|^2$ is
 an increasing (double) limit of bounded weakly continuous functions it follows that \eqref{e2.15}
 holds for  $\widetilde{\mu_t}$ in place of $\widetilde{\mu^\alpha_t}$ for each $t\in[s,T]\cap\Q.$
 Hence \cite[Theorem 8.6.7]{B2007} also applies to this family in $\mathcal P(H).$ In particular, for each
$t\in[s,T]\setminus\Q$ there exist $r_n(t)\in[s,T]\cap\Q,\;n\in \N,$ converging to $t$
and $\mu_t\in \mathcal P(H)$ such that $ \widetilde{\mu}_{r_n(t)}\to \mu_t$ $\tau_w$-weakly
as $n\to \infty$. We claim:
\begin{equation}
\label{e2.18}
\mu_t^{\alpha_n}\to\mu_t\;\mbox{\rm $\tau_w$-weakly as $n\to \infty$}\ \forall\;t\in[s,T]\setminus\Q.
\end{equation}
So, fix $t\in[s,T]\setminus\Q$ and suppose that $\{\mu_t^{\alpha_n}\}$ does not weakly converge
to $\mu_t$. Then by \eqref{e2.15} and \cite[Theorem 8.6.7]{B2007}  there exists a
subsequence $\{\alpha_{n_k}\}$ and $\nu\in \mathcal P(H)\setminus\{\mu_t\}$ such that
$\mu_t^{\alpha_{n_k}}\to\nu\;\tau_w$-weakly as $k\to \infty$.
Since $\mathcal E_A(H)$ is measure separating there exists
$\varphi\in \mathcal E_A(H)$ such that $\mu_t(\varphi)\neq \nu(\varphi)$.
On the other hand for all $n,k\in \N$ one has
\begin{multline*}
|\nu(\varphi)-\mu_t(\varphi)|
\le
|\nu(\varphi)-\mu^{\alpha_{n_k}}_t(\varphi)|
+\sup_{l\in\N}|\mu^{\alpha_{n_l}}_t(\varphi)-\mu^{\alpha_{n_l}}_{r_n(t)}(\varphi)|
\\
+
|\mu^{\alpha_{n_k}}_{r_n(t)}(\varphi)-\tilde\mu_{r_n(t)}(\varphi)|
+|\tilde\mu_{r_n(t)}(\varphi)-\mu_{t}(\varphi)|.
\end{multline*}
Since $\varphi$ is weakly continuous,
letting first $k\to\infty$ and then $n\to \infty$  it follows by \eqref{e2.17}
that $\mu_t(\varphi)=\nu(\varphi)$. This contradiction proves \eqref{e2.18}.
Letting $ \mu_t:=\tilde\mu_t$ for $t\in [s,T]\cap \Q$, by construction the first assertion
in Claim~1 follows for this family $\mu_t, t\in[s,T],$ in $\mathcal P(H)$.
Furthermore, \eqref{e2.15} then implies that
\begin{equation}
\label{e2.19}
\sup_{t\in[s,T]}\int_H|x|^2\mu_t(dx)<\infty.
\end{equation}
We have $\lim\limits_{n\to\infty}\mu^{\alpha_{n}}_t(\varphi)=\mu_t(\varphi)$ for all $t\in [s,T]$
and all $\varphi\in \mathcal E_A(H)$. Hence from \eqref{e2.17} the
second assertion in Claim 1 follows first for $\varphi\in \mathcal E_A(H)$,
but then by \eqref{e2.19} and Lebesgue's dominated convergence theorem,
this remains true for all $u\in D(L_0)$.
By a monotone class argument, the last   assertion in Claim 1
is  then an easy consequence.
Hence to complete the proof of Claim 1  it remains to prove \eqref{e2.17}.
So, fix $\varphi\in \mathcal E_A(H)$. Then by \eqref{e2.9}, \eqref{e2.4},
Hypothesis 2.3 and Lemma \ref{l2.4} for $s\le t_1\le t_2\le T$
\begin{multline}
\label{e2.20}
|\mu^\alpha_{t_2}-\mu^\alpha_{t_1}
|\le \frac12\;\|\mbox{\rm Tr}\;[CD^2\varphi]\|_\infty|t_2-t_1|\\
+|t_2-t_1|^{1/2}\|AD\varphi\|_\infty
\left(\int_{t_1}^{t_2}\int_H|x|^2\mu^\alpha_{s'}(dx)ds'\right) ^{1/2}
\\
+|t_2-t_1|^{1/2}K\|D\varphi\|_\infty
\left(\int_{t_1}^{t_2}\int_H V^2(s',x)\zeta(dx)ds'\right) ^{1/2},
\end{multline}
where $\|\cdot\|_\infty$ denotes sup-norm on $H$.
Since obviously all three sup-norms in \eqref{e2.20} are finite,
\eqref{e2.17} now follows from \eqref{e2.12} and \eqref{e2.15}. $\Box$\medskip

\noindent{\bf Proof of Claim 2}.
 For $\delta\in (0,\frac12)$ as in Hypothesis \ref{h2.1}(iii)
 from \eqref{e2.14} and \eqref{e2.14'} with $t=T$ we obtain for some $C>0$
\begin{multline*}
\int_s^T|(-A)^\delta Y_\alpha(t)|^2dt
\\
\le C \|(-A)^{-1/2+\delta}\|\left(|x|^2+\int_s^T|F_\alpha(s',X_\alpha(s'))|^2(x)ds'   \right).
\end{multline*}
Resubstituting according to \eqref{e2.13}, taking expectation and using \eqref{e2.1} we deduce that
\begin{multline*}
\int_s^T\E|(-A)^\delta X_\alpha(t,s,x)|^2dt
\\
\le 2C\|(-A)^{-1/2+\delta}\|\left(|x|^2
+\int_s^TP^\alpha_{s,s'}|F_\alpha(s',X_\alpha(s'))|^2(x)ds'   \right) +2c_\delta T.
\end{multline*}
Hence using  \eqref{e2.4}, Hypothesis \ref{h2.2} and Lemma \ref{l2.4} we find
\begin{multline}
\label{e2.23'}
\int_s^T\int_H|(-A)^\delta x|^2\mu_t^\alpha(dx)dt
\\
\le C_1\left(\int_H|x|^2\zeta(dx)+\int_s^T\int_HV(s',x)\zeta(dx)ds'   \right)
\end{multline}
for some constant $C_1$ independent of $\alpha$. Since $(-A)^{-\delta }$ is compact, $(-A)^\delta $ has compact level sets in $H$,
hence by Prohorov's theorem the sequence of measures $\mu_t^{\alpha_n}(dx)dt$
with $\alpha_n$ from Claim 1 has a subsequence weakly convergent
to a finite measure $\mu(dt,dx)$ (of total mass $T$) on $[0,T]\times H$.
 For simplicity we denote this subsequence again by $\{\mu_t^{\alpha_n}(dx)dt\}$.
 But for  $\varphi\in\mathcal E_A(H)$ and $f\in C_b([0,T];\R)$ we have
\begin{multline*}
\ds\int_s^T\int_Hf(t)\varphi(x)\mu_t(dx)dt
\ds=\int_s^Tf(t)\lim_{n\to\infty}\int_H\varphi(x)\mu_t^{\alpha_n}(dx)dt\\
\ds=\lim_{n\to\infty}\int_s^T\int_Hf(t)\varphi(x)\mu_t^{\alpha_n}(dx)dt
\ds=\int_s^T\int_Hf(t)\varphi(x)\mu_t(dx)dt,
\end{multline*}
where we used the weak continuity of $\varphi$ and
Lebesgue's dominated convergence theorem. From this it follows that $\mu(dt,dx)=\mu_t(dx)dt$.
The last part of Claim 2, i.e. \eqref{e2.12'},  follows from Lemma \ref{l2.4},
\eqref{e2.15}, \eqref{e2.23'} and the lower semicontinuity of $V+|(-A)^\delta \cdot |^2+|\cdot|^2$.
The proof of Claim 2 is complete. $\Box$\medskip

\noindent{\bf Proof of Claim 3}. We first note that
\eqref{e1.5'} is already verified because of \eqref{e2.11}. Furthermore, every $h\in C^1([0,T];D(A))$
can be written as a uniform limit of piecewise affine $h_n\in C([0,T];D(A)),\;n\in\N$,
uniformly bounded by $\|h'\|_\infty T$, e.g. by simply writing
$$
h(t)=h(0)+\int_0^th'(s)ds
$$
and approximating the integral by Riemannian sums. It then follows by Remark \ref{r1.1}(iii) and \eqref{e2.12'} by approximation and linearity that $\mu_t(dx)$ satisfies the Fokker--Planck equation \eqref{e1.5} or equivalently \eqref{(1.8)} if and only if it does so for all $u\in D(L_0)$ such that
$$
u(t,x)=\phi(t)e^{i\langle h(t),x  \rangle},\quad x\in H,\; t\in[0,T],
$$
with $\phi\in C^1([0,T];\R)$ and piecewise affine $h\in C([0,T];D(A))$. So, let us fix   such a
function $u\in D(L_0)$. Since \eqref{e1.5} and \eqref{(1.8)} are equivalent we know by \eqref{e2.9} that for all $n\in \N$
$$
\int_s^T\int_HL_{\alpha_n}u(t,x)\mu_t^{\alpha_n}(dx)dt=-\int_s^Tu(s,x)\zeta(dx)
$$
with $\alpha_n$ as in Claims 1,~2. Therefore, by Claim 2,
to show that \eqref{(1.6b)} holds
for $\mu_t(dx)dt$  it suffices to prove that for all $g\in C_b([s,T]\times H)$
\begin{equation}
\label{e2.21}
\lim_{n\to\infty}\int_s^T\int_HF^h_{\alpha_n}(t,x)g(t,x)\mu_t^{\alpha_n}(dx)dt=
\int_s^T\int_HF^h (t,x)g(t,x)\mu_t(dx)dt,
\end{equation}
where
$$
F^h_\alpha(t,x):=   \langle h(t),F_\alpha(t,x)
\rangle+\frac{\langle Ah(t),x  \rangle}{1+\alpha|\langle Ah(t),x  \rangle|},
$$
$$
F^h(t,x)\colon \,
= \langle h(t),F(t,x)   \rangle+\langle Ah(t),x  \rangle.
$$
 We note that $F^h_\alpha$ is continuous on $[s,T]\times H$ because of \eqref{e2.11}
 and because $ h$ is piecewise affine.  For $\delta\in(0,1]$ we have

\begin{multline}\label{e2.22}
 \left|\int_s^T\int_HF^h_{\alpha_n}(t,x)g(t,x)\mu_t^{\alpha_n}(dx)dt
 - \int_s^T\int_HF^h (t,x)g(t,x)\mu_t(dx)dt \right|
\\
 \le \|g\|_\infty\int_s^T\int_H|F^h_{\alpha_n}(t,x)-F^h(t,x)|\mu_t^{\alpha_n}(dx)dt
\\
+\|g\|_\infty\int_s^T\int_H|F^h(t,x)-F_\delta^h(t,x)|
\mu_t^{\alpha_n}(dx)dt
\\
+
\|g\|_\infty\int_s^T\int_H|F^h(t,x)-F_\delta^h(t,x)|
\mu_t(dx)dt
\\
+
\left|\int_s^T\int_HF_\delta^h(t,x)g(t,x)\mu_t^{\alpha_n}(dx)dt- \int_s^T\int_HF_\delta^h(t,x)g(t,x)\mu_t(dx)dt \right|.
\end{multline}
By \eqref{e2.5} and the inequality
$$
\left|\frac{a}{1+\delta|a|}-a\right|\le \delta|a|^2,\quad\forall\; a\in \R,
$$
we can find $\gamma(h)>0$ such that for all $n\in\N$ and all $\alpha,\beta\in (0,1]$
\begin{equation}
\label{e2.23}
\begin{array}{l}
\ds\int_s^T\int_H|F^h_{\beta}(t,x)-F^h(t,x)|\mu_t^{\alpha}(dx)dt\\
\\
\ds\le \beta\gamma(h)
\int_s^T\int_H(|F(t,x)|^2+|x|^2)\mu_t^{\alpha}(dx)dt.
\end{array}
\end{equation}
By Hypothesis \ref{h2.2}, Lemma \ref{l2.4}, \eqref{e2.15} and \eqref{e2.12} the integral on the right hand side of \eqref{e2.23} is bounded by a constant independent of $\alpha$.
So,  letting  $n\to\infty$ and $\delta\to 0$ the first two terms in
\eqref{e2.22} converge to zero. Using the last part of Claim 2, by the same arguments we deduce that this also holds for the third term. The last summand
 on the right hand side of \eqref{e2.22} can be estimated for every $\delta\in(0,1]$ by
 \begin{equation}
\label{e2.24}
\begin{array}{l}
\ds\left|\int_s^T\int_HF_\delta^h(t,x)g(t,x)(\mu_t^{\alpha_n}(dx)-\mu_t(dx))dt\right|\\
\\
\ds+\delta\|g\|_\infty\gamma(h) \int_s^T\int_H|F(t,x)|^2(\mu_t^{\alpha_n}(dx)+\mu_t(dx))dt  .
\end{array}
\end{equation}
Since, as pointed out above, $F^h$ is continuous on $[0,T]\times H$, the first summand in \eqref{e2.24}
converges to zero as $n\to \infty$ by Claim 2. Arguing as before we see that the the second summand
is bounded by $\delta$ times a constant independent of $n$.
So, letting first $n\to \infty$ and then $\delta\to 0$,
also the last term on  the right hand side of \eqref{e2.22} converges to zero
 and thus \eqref{e2.21} is proved, which completes the proof of Claim 3. $\Box$

\section{Uniqueness and Chapman--Kolmogorov\\ equations}

First let us recall the uniqueness result
from \cite{BDPR08c} on solutions of Fokker--Planck equations. This result
 is proved under certain assumptions on the coefficients $A,F$ and $C$ in equation \eqref{e1.1}
 which differ from those  in Section 2. We start with recalling them first.
\begin{Hypothesis}
\label{h3.1}
\begin{enumerate}

\item[]

\item[{\rm(i)}] There is $\omega \in \R$ such that $ \langle Ax,x \rangle \le \omega
|x|^{2},\; \forall\;x\in D(A).$

\item[{\rm(ii)}]  $C\in L(H)$   is symmetric, nonnegative and such that the linear operator
$$
Q^{(\alpha)}_t:=\int_0^ts^{-2\alpha}e^{sA}Ce^{sA^{*}}ds
$$
is of trace class for all $t>0$ and some $\alpha\in(0,\infty)$.

\item[{\rm(iii)}] Setting $Q_t:=\int_0^te^{sA}Ce^{sA^{*}}ds$,
one has  $e^{tA}(H)\subset Q_t^{1/2}(H)$ for all $t>0$
and there is    $\Lambda_t\in L(H)$ such that
$Q_t^{1/2}\Lambda_t=e^{tA}$ and
$$
\gamma_\lambda:=\int_0^{+\infty}e^{-\lambda t}\|\Lambda_t\|dt<+\infty.
$$
\end{enumerate}
\end{Hypothesis}
\begin{Hypothesis}
\label{h3.2}
There exists a family $\{\bar F(t,\cdot)\}_{t\in[0,T]}$ of $m$-quasi-dissipative maps
$$
\bar F(t,\cdot)\colon\
D(\bar F(t,\cdot))\subset H\;\to 2^H,\quad t\in[0,T],
$$
i.e., for each $t\in[0,T]$ the domain $D(\bar F(t,\cdot))$ belongs to
 $\mathcal B(H)$ and there exists $K>0$ independent of $t$ such that
$$
\langle u-v,x-y   \rangle\le K|x-y|^2,\quad\forall\;x,y\in D(\bar F(t,\cdot)),\; \;u\in \bar F(t,x),\;v\in \bar F(t,y)
$$
and for every $\lambda>K$ one has
$$
\mbox{\rm Range}\;(\lambda-\bar F(t,\cdot)):=\bigcup_{x\in D(\bar F(t,\cdot))}(\lambda x-\bar F(t,x))=H,
$$
\end{Hypothesis}
such that for every $t\in[0,T]$ we have $D(F(t,\cdot))=D(\bar F(t,\cdot))$
 and for all $x\in D(\bar F(t,\cdot)$ one has
\begin{equation}
\label{e3.3}
F(t,x)\in \bar F(t,x)\quad \mbox{\rm and}\quad |F(t,x)|=\min_{y\in\bar F(t,x)}\,|y|.
\end{equation}
\begin{Remark}
\label{r3.3}
\em We recall that for $\bar F(t,\cdot)$ as above the set $\bar F(t,\cdot)$ is convex and closed, so that the minimum in \eqref{e3.3} exists and   is unique.
\end{Remark}\medskip

To formulate the uniqueness result from \cite{BDPR08c}, let us introduce, for
every $\zeta\in\mathcal P(H)$ and $s\in[0,T]$, the set $\mathcal M_{s,\zeta}$ of all finite measures $\nu$ on $[s,T]\times H$ which have the following properties:\medskip

\noindent (i) $\nu(dt,dx)=\nu_t(dx)dt$, where $\nu_t(dx),\;t\in[s,T],$ is
a kernel from $([s,T],$ $\mathcal B([s,T]))$ to $(H,\mathcal B(H))$,
$\nu_t\in \mathcal P(H)$ for every $t\in[s,T]$ and $\nu_t(D(F(t,\cdot))=1$
 for $dt$-a.e. $t\in[s,T]$;
 \medskip

\noindent (ii) $\ds\int_s^T\int_H(|x|^2+|F(t,x)|+|x|^2|F(t,x)|)\nu_t(dx)dt<\infty$;
\medskip

\noindent (iii) $\nu_t(dx)dt$ satisfies identity \eqref{e1.5} for all $u\in D(L_0)$.
\medskip

We note that (ii) above implies \eqref{e1.5'} so that $L_0u\in L^1([0,T]\times H,\nu)$ for
all $u\in D(L_0),\;\nu\in\mathcal M_{s,\zeta}$.\medskip

\begin{Theorem} {\rm (\cite[Theorem \ref{t3.6}]{BDPR08c})}
\label{t3.4}
Suppose Hypotheses \ref{h3.1} and \ref{h3.2}
are fulfilled and $\zeta\in\mathcal P(H),\;s\in[0,T]$. Then $\mathcal M_{s,\zeta}$ contains at most one element.
\end{Theorem}
\begin{Remark}
\label{r3.5}
\em  (i) As already mentioned in \cite{BDPR08c}, combining the above theorem with
the results in \cite{BDPR08b} (see, in particular,
\cite[Corollary 1]{BDPR08b}), under additional coercivity conditions on the drift
one obtains quite general existence and uniqueness results for the
Fokker--Planck equation \eqref{e1.5}, more precisely, that
$\mathcal M_{s,\zeta}$ contains exactly one element.
From this, in the same way as explained below, one can
obtain the Chapman--Kolmogorov equation \eqref{(1.9)} for the transition functions,
i.e. the solutions $p_{s,t}(x,dy)dt$ of \eqref{e1.5} for $\zeta=\delta_x$,
at least for Lebesgue's a.e. $(r,s)\in [0,T)\times [0,T),\;r<s$.
On the basis of \cite{BDPR08b}, however, we can only treat cases where Tr $C<\infty$ (unless one can enlarge the state space $H$ in an appropriate way). Using our results from Section 2 above, in the present paper we shall analyze the case  Tr $C=\infty,$ more precisely,  the case $C^{-1}\in L(H)$. \medskip

\noindent(ii) By \cite[Remark 2.25]{DP2004} we have that Hypothesis \ref{h2.1} implies Hypothesis \ref{h3.1}.
\end{Remark}

 \begin{Theorem}
\label{t3.6}
Let $s\in[0,T]$ and suppose that Hypotheses \ref{h2.1}, \ref{h2.1r},
\ref{h3.2}, and  \eqref{e2.11} are fulfilled. Furthermore, assume that
Hypothesis \ref{h2.2} is fulfilled with $V$ satisfying
\begin{equation}
\label{e3.4}
|x|^2\le V(t,x),\quad\forall\;  (t,x)\in D(F),\;t\ge s.
\end{equation}
Then, for every $\zeta\in\mathcal P(H)$ satisfying \eqref{e2.12},
the measure $\mu_t(dx)dt$ from Theorem \ref{t2.5} is the only element
in $\mathcal M_{s,\zeta}$. In particular, for each $t\in[s,T]$ we have
$$
\mu_t^\alpha\to \mu_t\quad \tau_w\mbox{\it -weakly as}\; \alpha\to 0
$$
and
$$
\mu_t^\alpha(dx)dt\to \mu_t(dx)dt\quad \mbox{\it weakly as}\; \alpha\to 0
$$
(rather than only for a subsequence), and also for all $\varphi\in \mathcal E_A(H)$
$$
\lim\limits_{\alpha\to 0}\sup_{t\in[s,T]}
|\mu_t^\alpha(\varphi)- \mu_t(\varphi)|=0,
 $$
in particular,
$$
\mu_t(\varphi)\to \int_H\varphi(x)\zeta(dx) \quad\mbox{\it as}\; t\to 0.
$$
\end{Theorem}
{\bf Proof}. By Remark \ref{r3.5}(ii) we can apply Theorem \ref{t2.5} to obtain a measure $\mu_t(dx)dt$ which, as stated there, satisfies the defining  properties (i) and (iii) of $\mathcal M_{s,\zeta}$. But also (ii) holds by \eqref{e2.12'} since by \eqref{e3.4}
$$
|x|^2\,|F(t,x)|\le V^2(t,x),\quad\forall\;  (t,x)\in[s,T]\times H.
$$
For the proof of the last part of the assertion, we first recall that
the family of measures
$\mu_t^\alpha(dx)dt,\;\alpha\in(0,1],$ is a weakly compact set
 of finite positive measures of mass $T$ by \eqref{e2.23'} and, for each $t\in[s,T]$,
  by (2.18)  and \cite[Theorem 8.6.7]{B2007},
  the family of measures
   $\mu_t^\alpha(dx),\alpha\in(0,1],$ is a $\tau_w$-weakly compact set in $\mathcal P(H)$.

In the proof of Theorem \ref{t2.5} it was shown that  every  sequence
converging to zero in $(0,1]$
has a subsequence $\{\alpha_n\}$ such that $\mu_t^{  \alpha_n},\;t\in[s,T],$
satisfy Claims 1-3. However, as shown above,
  their corresponding limits $\mu_t(dx)dt$ must all coincide as
  measures on $[s,T]\times H$. Since all these limits have the
  property that $t\mapsto \mu_t(\varphi)$ is continuous on $[s,T]$ for each
   $\varphi\in \mathcal E_A(H)$ and the latter set is measure separating,
   it follows that for all these limits also measures $\mu_t,\;t\in[s,T],$
   are uniquely determined. Hence the last parts of the assertion also follow. $\Box$\medskip

As we shall see in the last section the additional assumption \eqref{e3.4} above is satisfied in many cases.\medskip

Now we turn to the Chapman--Kolmogorov equations \eqref{(1.9)}.
Let the assumptions in Theorem \ref{t3.6} hold  for all $s\in[0,T]$
  with the same function $V$ and set
$$
H_0:=\left\{x\in H:\;\int_0^TV^2(t,x)dt<\infty\right\}
$$
Then $H_0\in\mathcal B(H)$ and by Theorem \ref{t2.5} for every $x\in H_0$
and $s\in[0,T]$ there exists a measure $p_{s,t}(dx)dt$ on $[0,T]\times H$
having the properties listed in Theorem \ref{t2.5} with $\zeta=\delta_x$,
in particular, solving the Fokker--Planck equation \eqref{e1.5} with this initial
 condition for all $t\in[s,T]$.

\begin{Lemma}
\label{l3.7}
Let the assumptions of Theorem \ref{t3.6} hold for all $s\in[0,T]$  with the same function $V$ and let $s\in(0,T]$. Then for every $f\in \mathcal B_b(H)$ the map
$$
(t,x)\mapsto\one_{H_0}(x)\int_Hf(y)p_{s,t}(x,dy),\quad t\in [s,T],\;x\in H,
$$
is $\mathcal B([s,T])\times \mathcal B(H)$-measurable and for each $\varphi\in \mathcal E_A(H)$
$$
\lim_{t\to 0}\int_H\varphi(y)p_{s,t}(x,dy)=\varphi(x),\quad\forall\; x\in H_0.
$$
\end{Lemma}
{\bf Proof}. For all $\alpha\in(0,1]$ and $x\in H$ let $p^\alpha_{s,t}(x,dy)$ be the probability measure defined by $p^\alpha_{s,t}(x,A):=P^\alpha_{s,t}\one_A(x).$ Then for $\alpha_n:=\frac1n,\;n\in\N,$ it follows by the last part of Theorem \ref{t3.6} that for each $t\in [s,T],\;x\in H_0$ and $\varphi\in \mathcal E_A(H)$
$$
\int_H\varphi(y)p_{s,t}(x,dy)=\lim_{n\to\infty}
\int_H\varphi(y)p^{\alpha_n}_{s,t}(x,dy).
$$
Since the functions on the right are $\mathcal B([s,T])\times \mathcal B(H)$-measurable for each $n\in \N$ and $H_0\in \mathcal B(H)$, the first assertion is proved for $f=\varphi\in \mathcal E_A(H).$ For general $f\in \mathcal B_b(H)$ it then follows by a monotone class argument. The second assertion follows from the last part of Theorem \ref{t3.6}. $\Box$
 \begin{Theorem}
\label{t3.8}
Let the assumptions of Theorem \ref{t3.6} hold for all $s\in[0,T]$ with the same function $V$. Let $0\le r<t\le T$ and $p_{s,t}(x',dy),\;x'\in H_0,$ be as above. Then for every $x\in H_0$, $s\in(r,t)$ such that $p_{r,s}(x,H_0)=1$ we have
\begin{equation}
\label{e3.6}
\int_Hp_{s,t}(x',dy)p_{r,s}(x,dx')=p_{r,t}(x,dy),
\end{equation}
i.e. for all $f\in \mathcal B_b(H)$
$$
\int_H\int_Hf(y)p_{s,t}(x',dy)p_{r,s}(x,dx')=\int_Hf(y)p_{r,t}(x,dy),
$$
i.e. the Chapman--Kolmogorov equation holds.
\end{Theorem}

\begin{Remark}
\label{r3.9}
\em Let us discuss some conditions implying that
\begin{equation}
\label{e3.8}
p_{r,s}(x,H_0)=1,\quad\forall\; x\in H_0.
\end{equation}
Suppose that $D(F)=[0,T]\times Y$ for some set $Y\in \mathcal B(H)$. Then,
since $p_{r,s}(x,dy) ds\in \mathcal M_{r,\delta_x}$, we
know by its defining property (i) (stated before Theorem~3.4)
 that \eqref{e3.8} holds for $ds$-a.e. $s\in[r,T]$. To have it for all $s\in[r,T]$ let us assume that $V^2(\cdot,x)\in L^1(0,T;\R)$ for all $x\in Y$, hence $H_0=Y$, which is e.g. the case in our applications in Section 4 below.  We then know by \eqref{e2.10} that
\begin{equation}
\label{e3.9}
P^\alpha_{r,s}V^2(s,\cdot)(x)\le KV^2(s,x)<\infty, \quad\forall\;
x\in H_0,\;s\in[r,T],\;\alpha\in(0,1].
\end{equation}
Fix $x\in H_0$, $s\in[r,T]$. By construction (see the proof of Theorem \ref{t2.5}),
for any sequence $\alpha_n\to 0$, we know that
\begin{equation}
\label{e3.10}
\lim_{n\to \infty}p^{\alpha_n}_{r,s}(x,\cdot)=
p^{\alpha}_{r,s}(x,\cdot)\quad \tau_w\mbox{\rm -weakly},
\end{equation}
where $p^{\alpha_n}_{r,s}$ are as defined in the proof of Lemma \ref{l3.7}. Then we have
\begin{enumerate}
\item[(a)] If $V^2(s,\cdot)$ is an increasing limit of a sequence of weakly continuous functions (which is e.g. the case in our applications in Section 4 below), then
\eqref{e3.8} holds.

\item[(b)] If $V^2(s,\cdot)$ has compact level sets in the norm topology of $H$, then
\eqref{e3.8} holds.
\end{enumerate}

Property (a) follows immediately from \eqref{e3.9} and \eqref{e3.10}
 since we have
 $H\setminus H_0\subset\{V(s,\cdot)=\infty\}.$
 In case (b) one only has to note that by Prohorov's theorem it follows that the sequence of
 measures
 $p^{\alpha_n}_{r,s}(x,\cdot),\;n\in \N,$ is relatively compact also in
 the weak topology, so by \eqref{e3.10} it is even weakly convergent to $p_{r,s}(x,\cdot)$.

Since by the assumption in (b) the function $V^2(t,\cdot)$ is lower semicontinuous on $H$,
 \eqref{e3.9} implies by letting $n\to\infty$ that
$$
\int_HV^2(s,y)p_{r,s}(x,dy)\le KV^2(s,x),\quad\forall\;  x\in H_0,\;s\in[r,T].
$$
 Hence \eqref{e3.8} follows since $H\setminus H_0\subset\{V(s,\cdot)=\infty\}.$
\end{Remark}\medskip

\noindent{\bf Proof of Theorem \ref{t3.8}}. Let $x\in H_0$ and $u\in D(L_0)$.
Then for all  points $x'\in H_0$ and $t\in [s,T]$ one has
$$
\int_H u(t,y)p_{s,t}(x',dy) =u(s,x')+\int_s^t\int_HL_0u(s',y)p_{s,s'}(x',dy)ds'.
$$
Integrating with respect to $p_{r,s}(x,dx')$ and using Fubini's theorem (which is justified by Lemma \ref{l3.7})
we obtain for all $\;t\in[s,T]$ that
\begin{multline*}
\int_H u(t,y)\int_Hp_{s,t}(x',dy)p_{r,s}(x,dx')\\
 =\int_Hu(s,x')p_{r,s}(x,dx')+\int_s^t\int_HL_0u(s',y)\int_Hp_{s,s'}(x',dy)p_{r,s}(x,dx')ds'\\
 =u(r,x)+\int_r^s\int_HL_0u(s',y)p_{r,s'}(x,dx')ds'\\+\int_s^t\int_HL_0u(s',y)\int_Hp_{s,s'}(x',dy)p_{r,s}(x,dx')ds',
\end{multline*}
where we used that $p_{r,s}(x,dx')$ solves \eqref{e1.5} in the last
 equality. Hence defining the measures
$$
\mu^{(s)}_{r,s'}(x,dy):=\one_{[r,s]}(s')p_{r,s'}(x,dy)
+\one_{(s,T]}(s')\int_Hp_{s,s'}(x',dy)p_{r,s}(x,dx'),
$$
where $s'\in[s,T]$,
we have by the last part of Lemma \ref{l3.7} that for all $\varphi\in\mathcal E_A(H)$
the function
$$
s'\mapsto\int_H\varphi(y)\mu^{(s)}_{r,s'}(x,dy)
$$
is continuous on $[r,T]$ and $\mu^{(s)}_{r,s'}(x,dy)ds'$ satisfies the
Fokker--Planck equation \eqref{e1.5}, with $r,\;\delta_x$ in place of
 $s,\zeta$, respectively, and enjoys the defining properties (i) and (ii)
 for $\mathcal M_{r,\delta_x}$. But as noted above,
 we also have $p_{r,s'}(x,dy)ds'\in \mathcal M_{r,\delta_x}$, hence
$$
\mu^{(s)}_{r,s'}(x,dy)=p_{r,s'}(x,dy),\quad\forall\;  s'\in [r,T].
$$
In particular, \eqref{e3.6} holds. $\Box$

\section{Applications}
Let $H=L^2(0,1):=L^2((0,1),d\xi)$ and let $A\colon\, D(A)\subset H\to H$ be defined by
$$
Ax(\xi)=\partial^2_\xi x(\xi),\;\xi\in (0,1),\quad D(A)=H^2(0,1)\cap H^1_0(0,1),
$$
where $\partial_\xi=\frac{d}{d\xi},$  $\partial^2_\xi=\frac{d^2}{d\xi^2}.$

We would like to mention here that what is done below generalizes to the case where $(0,1)$ is replaced by an open set $\mathcal O$ in $\R^d,\;d\ge 1.$ One has only to replace the operator $C$ below by $A^{-\delta}$ with properly chosen $\delta>0,$
depending on the dimension $d$.

Let $D(F):=[0,T]\times L^{2m}(0,1)$ and for $(t,\xi)\in D(F)$
$$
F(t,x)(\xi):=f(\xi,t,x(\xi))+h(\xi,t,x(\xi)),\quad\;\xi\in (0,1).
$$
Here $f,h:(0,1)\times[0,T]\times\R\to \R$ are   functions such that
 for every $\;\xi\in (0,1)$ the maps $f(\xi,\cdot,\cdot),\; h(\xi,\cdot,\cdot)$ are
 continuous on $(0,T)\times\R$ and have the following properties:\medskip

\begin{enumerate}
\item[(f1)] (``polynomial bound''). There exist   $m\in\N$ and
a nonnegative function $c_1\in L^2(0,T)$  such that for all $t\in(0,T)$,
 $z\in\R,\;\xi\in (0,1)$ one has
$$
|f(\xi,t,z)|\le c_1(t)(1+|z|^m),
$$
also assuming without loss of generality that $m$ is odd.

\item[(f2)] (``quasi-dissipativity'').  There is a nonnegative function
$c_2\in L^1(0,T)$  such that for all $t\in[0,T],$ $z_1,z_2\in\R,\;\xi\in (0,1)$ one has
$$
(f(\xi,t,z_2)-f(\xi,t,z_1))(z_2-z_1)\le c_2(t)|z_2-z_1|^2.
$$

\item[(h1)] (``linear growth'').  There exists
a nonnegative function $c_3\in L^2(0,T)$ such that
for all $t\in[0,T],$ $z\in\R,\;\xi\in (0,1),$ one has
$$
|h(\xi,t,z)|\le c_3(t)(1+|z|).
$$
\end{enumerate}
Finally, let $C\in L(H)$ be symmetric, nonnegative and such that $C^{-1}\in L(H)$.

It is worth noting that it is not known whether
under these assumptions  the stochastic differential equation \eqref{e1.1} has a solution.

We set $Y:=D(F)=L^{2m}(0,1)$ and prove that Hypotheses \ref{h2.1}-\ref{h2.2}
and condition \eqref{e2.11} are fulfilled.
So, we can apply Theorem \ref{t2.5} to get existence of solutions to
 the Fokker-Planck equation \eqref{e1.5} in this situation.

Note that  Hypothesis \ref{h2.1} holds with $\omega=-\pi^2$ because $A^{-1}$
is of trace class. Furthermore,
 for $\alpha\in(0,1]$ and $(t,x)\in[0,T]\times H$
we set
  \begin{equation}
\label{e4.1}
F_\alpha(t,x):=\frac{F(t,x)(\xi)}{1+\alpha|F(t,x)(\xi)|},\quad \xi\in (0,1).
\end{equation}
Then $F_\alpha$ has all properties mentioned in
 Hypothesis \ref{h2.1r} and \eqref{e2.11} also holds by Lebesgue's
 dominated convergence theorem.
Set
\begin{equation}
\label{e4.2}
V(t,x):=\left\{\begin{array}{l}
2(c_1(t)+c_3(t))(1+|x|^m_{L^{2m}(0,1)})\quad\mbox{\rm if}\;(t,x)\in D(F)=[0,T]\times Y,\\
+\infty\quad\mbox{\rm otherwise}.
\end{array} \right.
\end{equation}
We are going to prove that Hypothesis \ref{h2.2} is fulfilled
for this function $V$ for all $s\in[0,T]$. First observe, that by (f1) and (h1) one has
\begin{equation}
\label{e4.3}
|F(t,x)|\le V(t,x)<\infty\quad\forall\; (t,x)\in D(F)=[0,T]\times Y.
\end{equation}
Furthermore, \eqref{e2.10} follows for all $s\in[0,T)$ from the next proposition.

\begin{Proposition}
\label{p4.1}
Let $\alpha\in(0,1]$, $x\in Y=L^{2m}(0,1)$  and  $s\in[0,T)$.
Let $X_\alpha(t,s,x),\;t\in [s,T],$ be the martingale solution
 of the approximating stochastic
differential equation \eqref{e2.6} started at $x$ at time $s$. Then there exists $C>0$ such that
\begin{equation}
\label{e4.4}
\E\left(|X_\alpha(t,s,x)|^{2m}_{L^{2m}(0,1)}   \right)\le C\left(1+|x|^{2m}_{L^{2m}(0,1)}   \right), \quad\forall\;t\in[s,T].
\end{equation}
\end{Proposition}
{\bf Proof}. We first note that   by (f1), (f2), and (h1),
for all $y,z\in\R$, $t\in[0,T], \xi\in (0,1),$ one has
 \begin{equation}
\label{e4.5}
\begin{array}{l}
f(\xi,t,y+z)y+h(\xi,t,y+z)y\\
\\
=(f(\xi,t,y+z)-f(\xi,t,z))y+f(\xi,t,z)y+h(\xi,t,y+z)y\\
\\
\le c_2(t)|y|^2+c_1(t)(1+|z|^m)|y|+c_3(t)(1+|y|+|z|)|y|\\
\\
\le c(t)(1+|y|^2+|z|^m|y|),
\end{array}
\end{equation}
where
$$
c=c_1+c_2+2c_3\in L^1(0,T).
$$
Setting
\begin{equation}
\label{e4.6}
Y_\alpha(t):=X_\alpha(t,s,x)-W_A(s,t),\quad t\in[s,T],
\end{equation}
\eqref{e2.6} reduces to
$$
\left\{\begin{array}{l}
\ds\frac{d}{dt}\;Y_\alpha(t)=AY_\alpha(t)+F_\alpha(t,X_\alpha(t,s,x)),\quad t\in[s,T],\\
\\
Y_\alpha(s)=x.
\end{array} \right.
$$
Here the equation is again meant in the  mild sense.
Now multiplying both sides of the first equation by $(Y_\alpha(t))^{2m-1}$
we obtain (after integration by parts) that for $t\in[s,T]$ one has
\begin{multline*}
\frac1{2m}\;\frac{d}{dt}\;\int_\mathcal O|Y_\alpha(t)|^{2m}d\xi
+(2m-1)\int_\mathcal O|Y_\alpha(t)|^{2m-2}|\partial_\xi Y_\alpha(t)|^2d\xi
\\
=\int_\mathcal OF_\alpha(t,Y_\alpha(t)+W_A(s,t))Y_\alpha(t)^{2m-1}d\xi,
\end{multline*}
where $\mathcal O:=(0,1).$ Taking into account \eqref{e4.1} and \eqref{e4.5}
we deduce that for $t\in[s,T]$ one has
\begin{multline*}
\frac1{2m}\;\frac{d}{dt}\;\int_\mathcal O|Y_\alpha(t)|^{2m}d\xi
\\
\le c(t)\int_\mathcal O[1+|Y_\alpha(t)|^2
+|W_A(s,t)|^m|Y_\alpha(t)|]|Y_\alpha(t)|^{2m-2}d\xi
\\
\le c(t)\int_\mathcal O\left[1+\left( 1+\frac{2m-1}{2m}  \right)|Y_\alpha(t)|^{2m}+\frac1{2m}\;|W_A(s,t)|^{2m^2}\right]d\xi,
\end{multline*}
which implies that for $c_m:=4m$ and
$$
\kappa:=1+\sup_{(t,\xi)\in [0,T]\times(0,1)}|W_A(s,t)(\xi)|
$$
one has
$$
\frac{d}{dt}\; |Y_\alpha(t)|_{L^{2m}(0,1)}^{2m}\le c_mc(t)\left(\kappa^{2m^2}+ |Y_\alpha(t)|_{L^{2m}(0,1)}^{m}  \right) ,\quad t\in [s,T].
$$
Applying a variant of Gronwall's lemma we arrive at
\begin{multline*}
|Y_\alpha(t)|_{L^{2m}(0,1)}^{2m}\le
\exp\biggl(c_m\int_s^tc(r)dr\biggr)|x|_{L^{2m}(0,1)}^{2m}
\\
+\kappa^{2m^2}c_m\int_s^t \exp\biggl(c_m\int_r^tc(r')dr'\biggr) c(r)dr.
\end{multline*}
Hence
$$
|Y_\alpha(t)|_{L^{2m}(0,1)}^{2m}\le e^{c_m |c|_{L^1(0,T)}}\left(|x|_{L^{2m}(0,1)}^{2m}+\kappa^{2m^2}c_m |c|_{L^1(0,T)}\right).
$$
However, according to \cite[Theorem 4.8(iii)]{DP2004}, we have
$$
\gamma_M:=\E(\kappa^M)<\infty,
$$
hence resubstituting according to \eqref{e4.6} we obtain
\begin{multline*}
\E|X_\alpha(t,s,x)|_{L^{2m}(0,1)}^{2m}
\\
\le 2^{m-1}\gamma_{2m}
+2^{m-1}e^{c_m|c|_{L^1(0,T)}}\left(|x|_{L^{2m}(0,1)}^{2m}
+\gamma_{2m^2}c_m |c|_{L^1(0,T)}\right),
\end{multline*}
and \eqref{e4.4} follows. $\Box$

Since $c_1,c_3\in L^2(0,T)$, Theorem \ref{t2.5} now applies to all $\zeta\in\mathcal P(H)$ such that
\begin{equation}
\label{e4.8}
\int_H|x|_{L^{2m}(0,1)}^{2m}\zeta(dx)<\infty.
\end{equation}
Now let us turn to uniqueness and the Chapman--Kolmogorov equations. Let $h\equiv 0$ and
$c_2\equiv $ const. Then $f$ is quasi-dissipative; in fact,
each $F(t,\cdot)$ with domain $Y$ is $m$-dissipative.
Hence Hypothesis \ref{h3.2} is also fulfilled. Furthermore,
the function $V$ defined in \eqref{e4.2} satisfies \eqref{e3.4},
hence Theorem \ref{t3.6} applies to give us uniqueness for solutions
of the Fokker--Planck equation \eqref{e1.5} for
every initial condition $\zeta\in\mathcal P(H)$ satisfying \eqref{e4.8},
in particular, for all $\zeta=\delta_x,\;x\in Y=L^{2m}(0,1).$

Furthermore, since $|\cdot|_{L^{2m}(0,1)}^{2m}$ is an increasing limit of a sequence of nonnegative weakly continuous functions on $L^2(0,1)$, by Remark
\ref{r3.9}, case (a), it follows that \eqref{e3.8}
holds for all $0\le r\le s\le T$ with $H_0:=L^{2m}(0,1).$ Hence by Theorem \ref{t3.8} the Chapman--Kolmogorov equation \eqref{(1.9)} holds for all
$x\in L^{2m}(0,1).$

\end{document}